\documentclass[a4paper,10pt]{article}

\usepackage[a4paper,left=3cm,right=3cm,top=2.5cm,bottom=2.5cm]{geometry}
\usepackage{hyperref}
\usepackage{subfig}
\usepackage{pgfplots}
\usepackage{pgfplotstable}
\usepackage{mathtools}
\usepackage{multicol}
\usepackage{comment}
\usepackage{booktabs}
\pgfplotsset{compat=1.5}
\usepackage{amssymb}
\usepackage{url}
\usepackage{bm}

\usepackage{pdfsync}
\usepackage{float}
\usepackage{tabularx}
\usepackage{enumerate}
\usepackage{array}
\usepackage{xspace}
\usepackage{siunitx}
\usepackage{authblk}

\usepackage[draft,inline,marginclue]{fixme}
\usepackage{mathrsfs}
\usepackage{color, colortbl}
\usepackage{multirow}

\usepackage{amsthm}
\usepackage{amsmath}
\usepackage{stmaryrd}
\usepackage[ruled,vlined,linesnumbered]{algorithm2e}
\usepackage{graphicx}
\usepackage{booktabs}
\usepackage{cleveref}

\FXRegisterAuthor{mt}{amt}{MT} \FXRegisterAuthor{fr}{afr}{FR}
\FXRegisterAuthor{al}{aal}{AL}

\theoremstyle{remark}

\newcommand{\x}{\ensuremath{\mathbf{x}}}

\newcommand{\y}{\ensuremath{\mathbf{y}}}

\newcommand{\X}{\ensuremath{\mathbf{X}}}

\newcommand{\R}{\ensuremath{\mathbb{R}}}
\newcommand{\K}{\ensuremath{\mathbf{K}}}

\newcolumntype{C}[1]{>{\centering\arraybackslash}m{#1}}

\definecolor{Gray}{gray}{0.9}

\begin{document}

\title{Multi-fidelity data fusion for the approximation of scalar functions with
low intrinsic dimensionality through active subspaces}
\author[]{Francesco~Romor\footnote{francesco.romor@sissa.it}}
\author[]{Marco~Tezzele\footnote{marco.tezzele@sissa.it}}
\author[]{Gianluigi~Rozza\footnote{gianluigi.rozza@sissa.it}}

\affil{Mathematics Area, mathLab, SISSA, via Bonomea 265, I-34136 Trieste,
  Italy}

\maketitle

\begin{abstract}
  Gaussian processes are employed for non-parametric regression in a Bayesian
  setting. They generalize linear regression embedding the inputs in a latent
  manifold inside an infinite-dimensional reproducing kernel Hilbert space. We
  can augment the inputs with the observations of low-fidelity models in order
  to learn a more expressive latent manifold and thus increment the model's
  accuracy. This can be realized recursively with a chain of Gaussian processes
  with incrementally higher fidelity. We would like to extend these
  multi-fidelity model realizations to case studies affected by a
  high-dimensional input space but with low intrinsic dimensionality. In these
  cases physical supported or purely numerical low-order models are still
  affected by the curse of dimensionality when queried for responses. When the
  model's gradient information is provided, the presence of an active subspace
  can be exploited to design low-fidelity response surfaces and thus enable
  Gaussian process multi-fidelity regression, without the need to perform new
  simulations. This is particularly useful in the case of data scarcity. In this
  work we present a multi-fidelity approach involving active subspaces and we
  test it on two different high-dimensional benchmarks.
\end{abstract}

\tableofcontents
\listoffixmes

\section{Introduction}
\label{sec:intro}

Every day more and more complex simulations are made possible thanks to the
high-performance computing facilities spread and to the advancements in
computational methods. Still the study and the approximation of high-dimensional
functions of interest represent a problem due to the curse of dimensionality and
the data scarcity because of limited computation budgets.

Gaussian processes (GP)~\cite{williams2006gaussian} have been proved as a
versatile and powerful technique for regression (GPR), classification, inverse
problem resolution, and optimization, among others.  In the last years several
studies and extensions have been proposed in the context of non-parametric and
interpretable Bayesian models. For a review on Gaussian processes and kernel
methods we suggest~\cite{kanagawa2018gaussian}, while for approximation methods
in the framework of GP see~\cite{quinonero2007approximation}. Progress has also
been made to deal with big data and address some memory limitations of GPR, as
for example by sparsifying the spectral representation of the
GP~\cite{lazaro2010sparse} or introducing stochastic variational inference for
GP models~\cite{liu2020gaussian}.

Multi-fidelity modelling has been proven effective in a heterogeneous set of
applications~\cite{kennedy2000predicting, forrester2007multi,
raissi2017inferring, bonfiglio2018multi, bonfiglio2018improving,
kramer2019multifidelity}, where expensive but accurate high-fidelity
measurements are coupled with cheaper to compute and less accurate low-fidelity
data. Recent advancements have been made for nonlinear autoregressive
multi-fidelity Gaussian process regression (NARGP) as proposed
in~\cite{perdikaris2017nonlinear}, and with physics informed neural networks
(PINNs)~\cite{raissi2019physics} in the context of multi-fidelity approximation
in~\cite{meng2020composite}.

The increased expressiveness of these models is achieved thanks to some kind of
nonlinear approach that extend GP models to non-GP processes with the
disadvantage of an additional computational cost. In this direction are focused
the following works which aim to obtain computationally efficient
heteroscedastic GP models with a variational inference
approach~\cite{lazaro2011variational} or employing a nonlinear
transformation~\cite{snelson2004warped}. This approach is extended to
multi-fidelity models departing from the linear formulation of Kennedy and
O'Hagan~\cite{kennedy2000predicting} towards deep Gaussian
processes~\cite{damianou2013deep} and NARPG.

When the models depend on a high-dimensional input space even the low-fidelity
approximations supported by a physical interpretation or a purely numerical
model order reduction suffer from the curse of dimensionality especially for the
design of high-dimensional GP models. Active subspaces
(AS)~\cite{constantine2015active, zahm2020gradient} can be used to build a
surrogate low-fidelity model with reduced input space taking advantage of the
correlations of the model's gradients when available. Reduction in parameter
space through AS has been proven successful in a diverse range of applications
such as: shape optimization~\cite{lukaczyk2014active, demo2020asga}, hydrologic
models~\cite{jefferson2015active}, naval and nautical
engineering~\cite{tezzele2018dimension, demo2018isope, mola2019marine,
tezzele2018model, tezzele2019marine}, coupled with intrusive reduced order
methods in cardiovascular studies~\cite{tezzele2018combined}, in CFD problems in
a data-driven setting~\cite{demo2019cras, tezzele2020enhancing}. A kernel-based
extension of AS for both scalar and vectorial functions can be found
in~\cite{romor2020kas}.

The aim of the present contribution is to propose a multi-fidelity regression
model which exploits the intrinsic dimensionality of high-dimensional functions
of interest and the presence of an active subspace to reduce the approximation
error of high-fidelity response surfaces. Our approach employs the design of a
NARGP using an AS response surface as low-fidelity model. In the literature the
multi-fidelity approximation paradigm has been adopted in a different way to
search for an active subspace from given high- and low-fidelity
models~\cite{lam2020multifidelity}.

The outline of this work is the following: in Section~\ref{sec:multifidelity} we
present the general framework of Gaussian process regression and in particular
the NARGP multi-fidelity approach; in Section~\ref{sec:as} we briefly present
the Active subspaces property we are going to exploit as low-fidelity model,
together with some error estimates for the construction of response surfaces;
Section~\ref{sec:mfas} is devoted to present how data fusion with Active
subspaces is performed with the aid of algorithms; in Section~\ref{sec:results}
we apply the proposed approach to the piston and Ebola benchmark models showing
the better performance achieved by the multi-fidelity regression; finally in
Section~\ref{sec:conclusions} we present our conclusions and we draw some future
perspectives.

\section{Multi-fidelity Gaussian process regression}
\label{sec:multifidelity}

In the next subsections we are going to present the Gaussian process regression
(GPR)~\cite{williams2006gaussian} technique which is the building block of
multi-fidelity GPR, and the nonlinear autoregressive multi-fidelity Gaussian
process regression (NARGP)~\cite{perdikaris2017nonlinear} we are going to use in
this work. The numerical methods are presented for the general case of several
levels of fidelity.

\subsection{Gaussian process regression}

Gaussian process regression~\cite{williams2006gaussian} is a supervised
technique to approximate unknown functions given a finite set of input/output
pairs $\mathcal{S} = \{ x_i, y_i \}_{i=1}^N$. Let $f: \mathcal{X} \subset \R^D
\to \R$ be the scalar function of interest. The set $\mathcal{S}$ is generated
through $f$ with the following relation: $y_i = f(x_i)$, which are the
noise-free observations. To $f$ is assigned a prior with mean $m(\x)$ and
covariance function $k(\x, \x^\prime ; \theta)$, that is $f (\x) \sim
\mathcal{GP} (m(\x), k(\x, \x^\prime ; \theta) )$. The prior expresses our
beliefs about the function before looking at the observed values. From now on we
consider zero mean $\mathcal{GP}$, $m(\x) = \mathbf{0}$, and we define the
covariance matrix $\K_{i,  j} = k(x_i, x_j ; \theta)$, with $\K \in \R^{N \times
N}$. To use the Gaussian process to make prediction we still need to find the
optimal values of the elements of the hyper-parameters vector $\theta$. We
achieve this by maximizing the log marginal likelihood:
\begin{equation}
\log p(\y | \x, \theta) = - \frac{1}{2} \y^T \K^{-1} \y -\frac{1}{2}
\log |\K| - \frac{N}{2} \log 2 \pi.
\end{equation}
Let $\x_*$ be the test samples, and $\K_{N*} = k(\x, \x_* ; \theta)$ be the
matrix of the covariances evaluated at all pairs of training and test samples,
and in a similar fashion $\K_{*N} = k(\x_*, \x ; \theta)$, and $\K_{**} =
k(\x_*, \x_* ; \theta)$. By conditioning the joint Gaussian distribution on the
observed values we obtain the predictions $f_*$ by sampling the posterior
\begin{equation}
f_* | \x_*, \x, \y \sim \mathcal{N} (\K_{*N} \K^{-1} \y, \K_{**} -
\K_{*N} \K^{-1} \K_{N*} ).
\end{equation}

\subsection{Nonlinear multi-fidelity Gaussian process regression}

We adopt the nonlinear autoregressive multi-fidelity Gaussian process regression
(NARGP) scheme proposed in~\cite{perdikaris2017nonlinear}. It extends the
concepts present in~\cite{kennedy2000predicting, le2014recursive} to nonlinear
correlations between the different fidelities available.

We introduce $p$ levels of increasing fidelities and the corresponding sets of
input/output pairs $\mathcal{S}_q = \{ x_i^q, y_i^q \}_{i=1}^{N_q} \subset
\mathcal{X} \times\R \subset\R^n \times \R$ for $q \in \{1, \dots, p \}$, where
$y_i^q = f_q (x_i^q)$. With $p$ we indicate the highest fidelity. We also assume
the design sets to have a nested structure: $\mathcal{S}_p \subset
\mathcal{S}_{p-1} \subset \cdots \subset \mathcal{S}_1$.

The NARGP formulation considers the following autoregressive multi-fidelity
scheme:
\begin{equation}
  f_q (\x) = g_q (\x, f_{* q-1} (\x) ) ,
\end{equation}
where $f_{* q-1} (\x)$ is the GP posterior from the previous inference level
$q-1$, and to $g_q$ is assigned the following prior:
\begin{align}
  g_q &\sim \mathcal{GP} ( f_q | \mathbf{0} , k_q ; \theta_q ) \\
  k_q & = k_q^\rho (\x, \x^\prime ; \theta_q^\rho) \cdot k_q^f ( f_{*
        q-1} (\x), f_{* q-1} (\x^\prime) ; \theta_q^f ) + k_q^\delta
        (\x, \x^\prime ; \theta_q^\delta),
\end{align}
with $k_q^\rho$, $k_q^f$, and $k_q^\delta$ squared exponential kernel functions.
With this scheme, through $g_q$ we can infer the high-fidelity response by
projecting the lower fidelity posterior to a latent manifold of dimension $D+1$.
This structure allows for nonlinear and more general cross-correlations between
subsequent fidelities.

A part from the first level of fidelity $q=1$ the posterior probability
distribution given the previous fidelity models is no longer Gaussian since the
inputs are couples $((\x, \x_{*}),(y_{q-1}(\x), f_{*q-1}(\x_{*}))$ where
$f_{q-1}$ is a Gaussian process $g_q \sim \mathcal{GP} ( f_q | \mathbf{0} , k_q;
\theta_q )$, the training set is $(\x,f_{q-1}(\x))$ and $\x_{*}$ is the new
input. So in order to evaluate the predictive mean and variance for a new input
$\x_{*}$ we have to integrate the usual Gaussian posterior $p(f_{*q}(\x_{*},
f_{*q-1}(\x_{\star}))|f_{*q-1}, \x_{*}, \x_{q}, \y_{q})$ explicited as
\begin{align}
&f_{*q}(\x_{*}, f_{*q-1}(\x_{\star}))|f_{*q-1}, \x_{*}, \x_{q}, \y_{q}\sim\mathcal{N} (\K^{q}_{*N} (\K^{q})^{-1} y_{q}, \K^{q}_{**} -
\K^{q}_{*N} (\K^{q})^{-1} \K^{q}_{N*} ),\\
&K_{*N}^{q}=k_{q}((\x_{*},f_{*q-1}(\x_{*}), (\x_{q-1}, y_{q-1});\theta),\\
&K_{N*}^{q}=k_{q}((\x_{q-1}, y_{q-1}), (\x_{*},f_{*q-1}(\x_{*});\theta),\\
&K^{q}=k_{q}((\x_{q-1}, y_{q-1}), (\x_{q-1}, y_{q-1});\theta),
\end{align}
over the Gaussian distribution of the prediction at the previous level
$f_{*q-1}(\x_{*})\sim\mathcal{N}(\mu_{x_{*}},\sigma_{\x_{*}})$. In practice the
following integral is approximated with recursive Monte Carlo in each fidelity
level
\begin{equation}
p(f_{*q}(\x_{*}, f_{*q-1}(\x_*)))=\int_{\mathcal{X}} p(f_{*q}(\x_{*},
f_{*q-1}(\x_*))|f_{*q-1}, \x_{*}, \x_{q}, \y_{q}) p(f_{*q-1}(\x_{*}))
d\x_{*}.
\end{equation}

\section{Active subspaces property}
\label{sec:as}

In this section we follow the work of P.G.
Constantine~\cite{constantine2015active}. Our aim is testing multi-fidelity
Gaussian process regression models to approximate objective functions which
depend on inputs sampled from a high-dimensional space. Low-fidelity models
relying on a physical supported or numerical model reduction (for example a
coarse discretization or a more specific numerical model order reduction) still
suffer from the high dimensionality of the input space.

In our approach we try to tackle these problematics by searching for a surrogate
(low-fidelity) model accounting for the complex correlations among the input
parameters that concur to the output of interest. With this purpose in mind we
query the model for an active subspace and employing the found active subspace
we design a response surface with Gaussian process regression.

\subsection{Response surface design with active subspaces}
Let us suppose the inputs are represented by an absolutely continuous random
variable $\X$ with probability density $\rho$ such that
$\text{supp}(\rho)=\mathcal{X}\subset\mathbb{R}^{m}$ where $m$ is the dimension
of the input space. If our numerical simulations provide also the gradients of
the samples for which the model is inquired for, we can approximate the
correlation matrix of the gradient with simple Monte Carlo as
\begin{equation}
\mathbb{E}_{\rho}[\nabla_{\x} f (\nabla_{\x}
f)^{T}]\approx\frac{1}{N}\sum_{i=1}^{N} \nabla_{\x} f(\X_{i})
(\nabla_{\x} f(\X_{i}))^{T} ,
\end{equation}
where $N$ is the number of sampled inputs. We then search for the highest
spectral gap $\lambda_{r}-\lambda_{r+1}$ in the sequence of ordered eigenvalues
of the approximated correlation matrix
\begin{equation}
\lambda_{1}\geq\dots\geq\lambda_{r}\geq\lambda_{r+1}\geq \dots \geq \lambda_{m}.
\end{equation}
The active subspace is the eigenspace corresponding to the eigenvalues
$\lambda_{1},\dots,\lambda_{r}$ and we represent it with the matrix
$\hat{W}_{1}\in\mathcal{M}(m\times r)$ which columns are the first $r$ active
eigenvectors. Then the response surface $\mathcal{R}$ is built with a Gaussian
process regression over the training set composed by $N_{\text{train}}$
input-output pairs $\{W_{1}^{T}\x_{i}, y_{i}\}_{i=1}^{N_{\text{train}}}$.

An a priori bound on the error can be proved~\cite{constantine2015active}, but
the whole approximation procedure considers additional steps for the evaluation
of the optimal profile $\mathbb{E}_{\rho}[f|\hat{W}_{1}^{T}\X]$ and its
approximation with Monte Carlo
$\overline{\mathbb{E}_{\rho}[f|\hat{W}_{1}^{T}\X]}$,
\begin{equation}
    f(\X) \approx \mathbb{E}_{\rho}[f|\hat{W}_{1}^{T}\X]\approx
    \overline{\mathbb{E}_{\rho}[f|\hat{W}_{1}^{T}\X]}\approx
    \mathcal{R}(\hat{W}_1\X).
\end{equation}
The mean square regression error is bounded a-priori by,
\begin{align*}
    \mathbb{E}_{\rho}&\left[(f(\X)-\mathcal{R}(\hat{W}_{1}^{T}\X))^2\right] \\
    &\leq \ C_1(1+N^{-1/2})^2\left(
      \epsilon(\lambda_1+\dots+\lambda_r)^{1/2}+(\lambda_{r+1}+\dots +
      \lambda_m)^{1/2}\right)^2+C_2\delta,
\end{align*}
where $C_{1}$ and $C_{2}$ are constants, $\epsilon$ quantifies the error in the
approximation of the true active subspace $W_{1}$ with $\hat{W}_{1}$ obtained
from the Monte Carlo approximation and $C_{2}\delta$ is a uniform bound on
\begin{equation}
\mathbb{E}_{\rho}[\overline{(\mathbb{E}_{\rho}[f|\hat{W}_{1}^{T}\X]} -
\mathcal{R}(\hat{W}_{1}^{T}\X)^{2} | W_{2}^{T}\X]\leq
C_{2}\delta.
\end{equation}

\section{Multi-fidelity data fusion with active subspaces}
\label{sec:mfas}
Our study is based on the design of a nonlinear autoregressive multi-fidelity
Gaussian process regression (NARGP)~\cite{perdikaris2017nonlinear} whose
low-fidelity level is learnt from a response surface built through the active
subspaces methodology. In fact we suppose that the model in consideration has
indeed a high dimensional input space but its intrinsic dimensionality is
sufficiently lower. This is often the case as shown by the numerous industrial
applications~\cite{jefferson2015active, tezzele2018combined,
tezzele2018dimension, lukaczyk2014active}.

The whole procedure requires the knowledge of an input/output high-fidelity
training set $\{(\x^H_i, y^H_i)\}_{i=1}^{N_H}\subset\R^{m}\times\R$, completed
by the gradients $\{dy^H_i\}_{i=1}^{N_H}\subset\R^{m}$ needed for the active
subspace's presence inquiry and a low-fidelity input set
$\{\x^L_i\}_{i=1}^{N_L}\subset\R^{m}$. We represent with $N_H, N_L$ the number
of high-fidelity and low-fidelity training set samples, respectively. Differently from
the usual procedure the low-fidelity outputs $\{y^L_i\}_{i=1}^{N_L}$ are
predicted with the response surface built thanks to the knowledge of the active
subspace through the dataset $\{(\hat{W}_{1}\x^H_i, y^H_i)\}_{i=1}^{N_H}$. At
the same time the response surface is also queried for the predictions $\{y^{H,
\text{train}}_i\}_{i=1}^{N_H}$ at the high-fidelity inputs
$\{\x^H_i\}_{i=1}^{N_H}$ that will be used for the training of the
multi-fidelity model. Now all the ingredients for the same procedure described
in~\cite{perdikaris2017nonlinear} are ready: the multi-fidelity model is trained
at the low-fidelity level with $\{(x^L_i, y^L_i)\}_{i=1}^{N_L}$ and at the
high-fidelity level with $\{((x^H_i, y^{H, \text{train}}_i),
y^H_i)\}_{i=1}^{N_H}$.

We remark that in this case the same high-fidelity outputs
$\{y^H_i\}_{i=1}^{N_H}\subset\R$ are used for the response surface training and
the high-fidelity training of the multi-fidelity model. In fact the outputs
$\{y^{H, \text{train}}_i\}_{i=1}^{N_H}$ predicted with the response surface are
equal to $\{y^H_i\}_{i=1}^{N_H}\subset\R$ since the response surface is a
Gaussian process with no noise trained on the dataset $\{(\hat{W}_{1}\x^H_i,
y^H_i)\}_{i=1}^{N_H}$ and queried for the same inputs
$\{\hat{W}_{1}\x^H_i\}_{i=1}^{N_H}$ for the predictions, that is $\{y^{H,
\text{train}}_i\}_{i=1}^{N_H}$. This results in the training of the
high-fidelity level of the multi-fidelity model with the dataset $\{((x^H_i,
y^{H, \text{train}}_i), y^H_i)\}_{i=1}^{N_H}=\{((x^H_i, y^H_i),
y^H_i)\}_{i=1}^{N_H}$.

A second procedure is developed where part of the high-fidelity inputs is used
only to train the response surface such that in general
$\{y^H_i\}_{i=1}^{N_H}\neq\{y^{H,\text{train}}_i\}_{i=1}^{N_H}$. The main
difference is that the surrogate low-fidelity model is built independently from
the high-fidelity level of the multi-fidelity model.

We expect that with the multi-fidelity approach, thanks to the nonlinear
fidelity fusion realized by the method, not only the lower accuracy of the
low-fidelity model will be safeguarded against, but also a hint towards the
presence of an active subspace will be transferred from the low-fidelity to the
high-fidelity level. In fact the low-fidelity GP regression model is built from
the predictions obtained with the $r$-dimensional response surface which
expressiveness is guaranteed by the additional assumption that the model under
investigation has a $r$-dimensional active subspace. So a part from the lower
computational budget and the reduced accuracy, our low-fidelity model should
transfer to the high-fidelity level the knowledge of the presence of an active
subspace when learning correlations among the inputs $\{x^H_i \}_{i=1}^{N_H}$,
the response surfaces predictions $\{y^{H, \text{train}}_i\}_{i=1}^{N_H}$ and
the high-fidelity targets $\{y^H_i\}_{i=1}^{N_H}$. The overhead with respect to
the original procedure~\cite{perdikaris2017nonlinear} is the evaluation of the
active subspace from the high-fidelity inputs.

The procedure is synthetically reviewed through
Algorithm~\ref{algo:mfas_dependent} for the use of the same high-fidelity
samples in the training of the response surface and of the second fidelity level
of the multi-fidelity model, and Algorithm~\ref{algo:mfas_independent} for the
use of independent samples. The main difference in the two procedures is the set
of samples with which the active subspace is computed.

\begin{algorithm}[H]
    \caption{NARGP with the same high-fidelity samples for AS response surface design}
    \label{algo:mfas_dependent}

    \SetKwInOut{Input}{input}\SetKwInOut{Output}{output}

    %\Indm
    \Input{high-fidelity inputs, outputs, gradients triplets $\{(\x^{H}_{i},
        y^{H}_{i}, dy^{H}_{i})\}_{i=1}^{N_H}\subset \R^{m}\times \R\times
        \R^{m}$,\\
    low-fidelity inputs $\{\x^{L}_{i}\}_{i=1}^{N_L}\subset \R^{m}$,\\
    training dataset $\{(\x_{i}^{\text{test}},
        y_{i}^{\text{test}})\}_{i=1}^{N_{\text{test}}}$,\\
    } \BlankLine \Output{multi-fidelity model, $g_{M}=((f_{H}|x^{H}_{i}, y^{H,
                \text{train}}_{i}),\ (f_{L}|x^{L}_{i}))\sim
                (\mathcal{GP}(f_{H}|m_{H}, \sigma_{H}),
                \mathcal{GP}(f_{L}|m_{L}, \sigma_{L}))$}

    %\Indp
    \BlankLine Compute the active subspace $\hat{W}_{1}$ with the high-fidelity
    gradients $\{dy_{i}^{H}\}_{i=1}^{N_H}$,\\
    Build the one-dimensional response surface $\mathcal{R}(\hat{W}_{1}\X)$ with
    a GP regression from $\{(\hat{W}_{1}\x^{H}_{i}, y^{H}_{i})\}_{i=1}^{N_H}$,\\
    Predict the low-fidelity outputs $\{y^{L}_{i}\}_{i=1}^{N_L}$ at
    $\{\x^{L}_{i}\}_{i=1}^{N_L}$ and the training high-fidelity inputs $\{y^{H,
    \text{train}}_{i}\}_{i=1}^{N_H}$ at $\{\x^{H}_{i}\}_{i=1}^{N_H}$ with the
    response surface,\\
    Train the multi-fidelity model at the low-fidelity level $g_{L}$ with the
    training dataset $\{(x^{L}_{i}, y^{L}_{i})\}_{i=1}^{N_L}$,\\
    Train the multi-fidelity model at the high-fidelity level $g_{H}$ with the
    training dataset $\{((x^{H}_{i}, y^{H, \text{train}}_{i}),
    y^{H}_{i})\}_{i=1}^{N_H}$

    %\Indp
    \BlankLine
\end{algorithm}
\clearpage
\begin{algorithm}[H]
    \caption{NARGP with additional independent high-fidelity samples for AS
      response surface design}
    \label{algo:mfas_independent}

    \SetKwInOut{Input}{input}\SetKwInOut{Output}{output}

    \Input{high-fidelity inputs, outputs, gradients triplets $\{(\x^{H}_{i},
        y^{H}_{i}, dy^{H}_{i})\}_{i=1}^{N_H}\subset \R^{m}\times \R\times
        \R^{m}$,\\
    low-fidelity inputs $\{\x^{L}_{i}\}_{i=1}^{N_L}\subset \R^{m}$,\\
    training dataset $\{(\x_{i}^{\text{test}},
        y_{i}^{\text{test}})\}_{i=1}^{N_{\text{test}}}$,\\
    AS procedure's inputs, outputs, gradients triplets $\{(\x^{\text{AS}}_{i},
    y^{\text{AS}}_{i}, dy^{\text{AS}}_{i})\}_{i=1}^{N_{\text{AS}}}\subset
    \R^{m}\times \R\times \R^{m}$,\\
} \BlankLine \Output{multi-fidelity model, $g_{M}=((f_{H}|x^{H}_{i}, y^{H,
\text{train}}_{i}),\ (f_{L}|x^{L}_{i}))\sim (\mathcal{GP}(f_{H}|m_{H},
\sigma_{H}), \mathcal{GP}(f_{L}|m_{L}, \sigma_{L}))$}

    \BlankLine Compute the active subspace $\hat{W}_{1}$ with the active
    subspaces dataset's gradients $\{dy_{i}^{\text{AS}}\}_{i=1}^{N_H}$,\\
    Build the one-dimensional response surface $\mathcal{R}(\hat{W}_{1}\X)$ with
    a GP regression from $\{(\hat{W}_{1}\x^{\text{AS}}_{i},
    y^{\text{AS}}_{i})\}_{i=1}^{N_{\text{AS}}}$,\\
    Predict the low-fidelity outputs $\{y^{L}_{i}\}_{i=1}^{N_L}$ at
    $\{\x^{L}_{i}\}_{i=1}^{N_L}$ and the training high-fidelity inputs $\{y^{H,
    \text{train}}_{i}\}_{i=1}^{N_H}$ at $\{\x^{H}_{i}\}_{i=1}^{N_H}$ with the
    response surface,\\
    Train the multi-fidelity model at the low-fidelity level $g_{L}$ with the
    training dataset $\{(x^{L}_{i}, y^{L}_{i})\}_{i=1}^{N_L}$,\\
    Train the multi-fidelity model at the high-fidelity level $g_{H}$ with the
    training dataset $\{((x^{H}_{i}, y^{H, \text{train}}_{i}),
    y^{H}_{i})\}_{i=1}^{N_H}$

    % \Indp
    \BlankLine

\end{algorithm}

Some final remarks are due. As in~\cite{perdikaris2017nonlinear} we assume that
the observations $\{y^{q}_{i}\}$ are noiseless for each level of fidelity $q$.
We employ Radial Basis Function kernels (RBF) with Automatic Relevance
Determination (ARD). The hyperparameters tuning is achieved maximizing the
log-likelihood with the gradient descent optimizer L-BFGD.

\section{Numerical results}
\label{sec:results}
We consider two different benchmark test problems for which a multi-fidelity
model will be built. The first is a $8$-dimensional model for the spread of
Ebola\footnote{The Ebola dataset was taken from
\url{https://github.com/paulcon/as-data-sets}.} and the second is a
$7$-dimensional model to compute the time it takes a cylindrical piston to
complete a cycle\footnote{The piston dataset was taken from
\url{https://github.com/paulcon/active_subspaces}.}. The library employed to
implement the NARGP model is Emukit~\cite{emukit2019} while for the active
subspace's response surface design we have chosen ATHENA\footnote{Available at
\url{https://github.com/mathLab/ATHENA}.}~\cite{athena2020} for the active
subspaces presence study and GPy~\cite{gpy2014} for the GP regression.

These tests have already been analyzed for the presence of an active subspace
and they indeed present a low intrinsic dimensionality. For each model we show
the sufficient summary plot along the one-dimensional active subspace
found, and the correlation among the low-fidelity level and the high-fidelity level of the
multi-fidelity model. We also present a comparison of the prediction error with respect to a
low-fidelity model (LF) represented by a GP regression on the low-fidelity
input/output dataset, a high-fidelity model (HF) represented by a GP regression
on the high-fidelity input/output dataset, and the proposed multi-fidelity model
(MF). In each test case the number of low-fidelity samples is $200$ and an error
study over the number of high-fidelity samples used is undergone. We
apply both the algorithms presented in the previous section. In
particular for Algorithm~\ref{algo:mfas_dependent}, the number of samples
used to find the active subspace is one third of the total number of
high-fidelity samples rounded down.

\paragraph{The piston model}
The algebraic cylindrical piston model appeared as a test for statistical
screening in~\cite{ben2007modeling}, while applications of AS to this model can
be found in~\cite{constantine2017global}. The scalar target function of interest
is the time it takes the piston to complete a cycle, and its computation
involves a chain of nonlinear functions. This quantity depends on $7$ input
parameters uniformly distributed. The corresponding variation ranges are taken
from~\cite{constantine2017global}.  The 10000 test points are samples with
Latin-Hypercube sampling.

\begin{figure}[ht!]
  \centering
  \includegraphics[trim=0 10 0 10, clip, width=.48\textwidth]{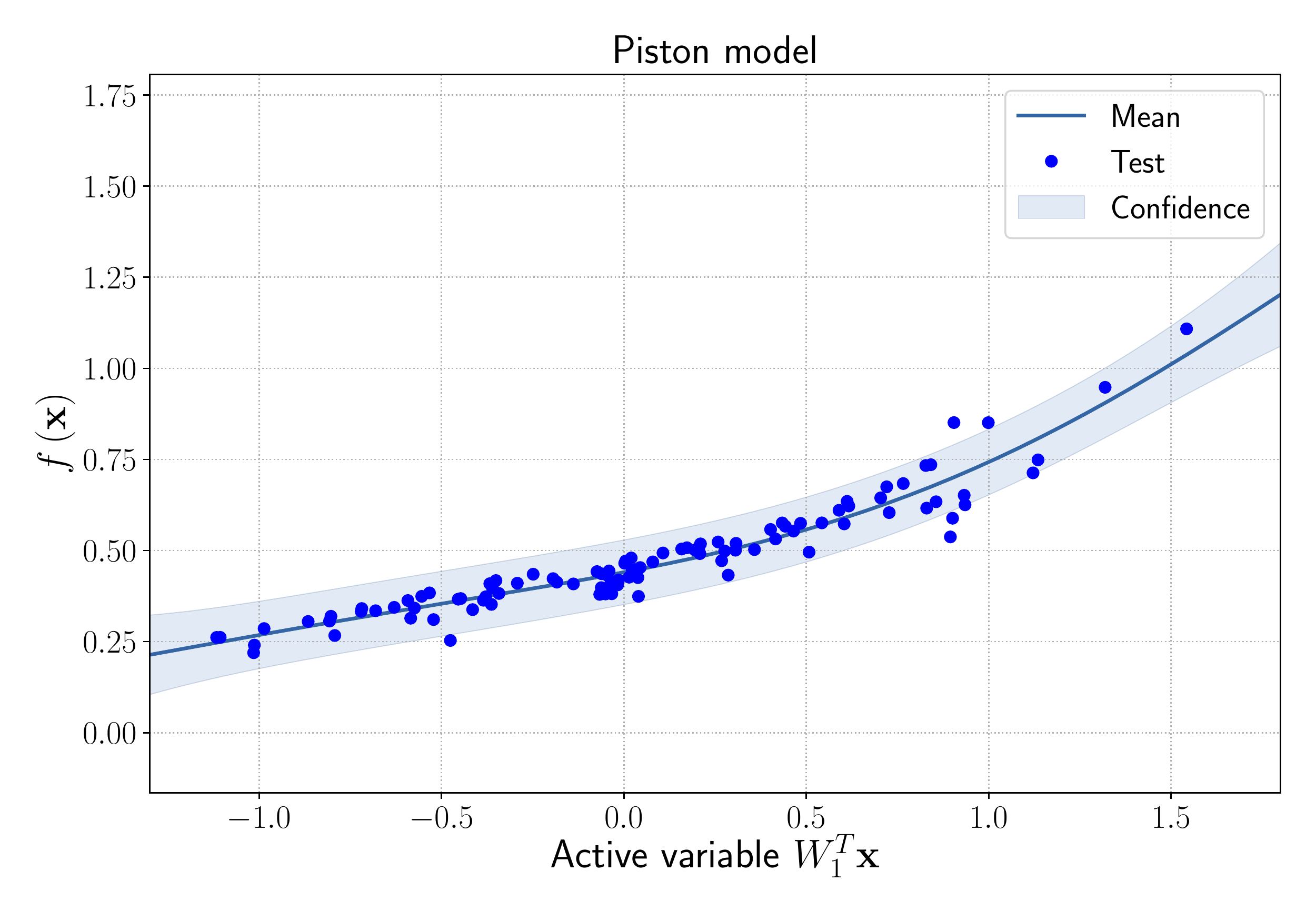}\hfill
  \includegraphics[trim=0 0 0 10, clip, width=.49\textwidth]{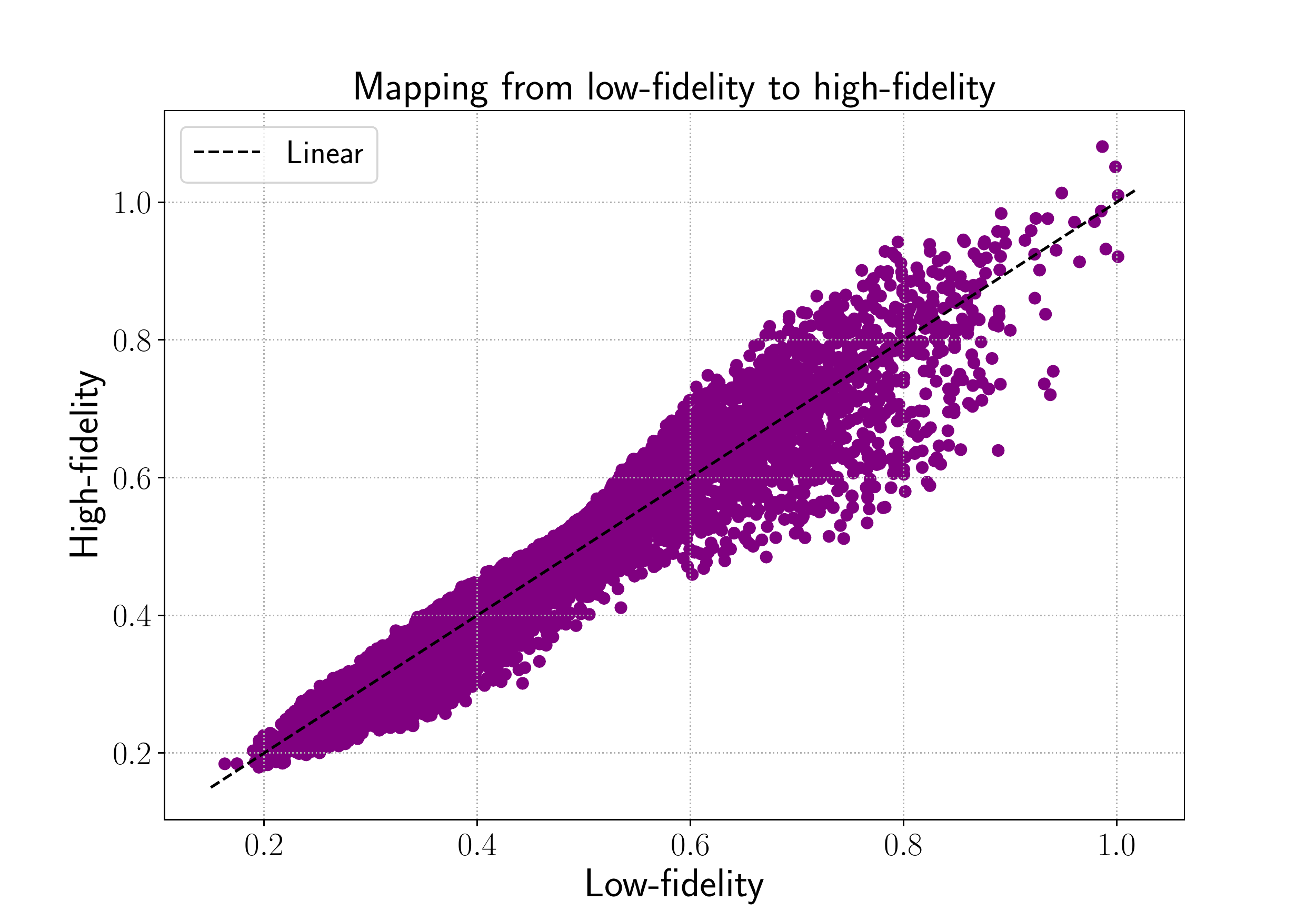}\\
  \caption{Left: sufficient summary plot of the surrogate model built
    with active subspaces. Right: correlation among the low-fidelity
    level and the high-fidelity level of the multi-fidelity model.}
  \label{fig:piston_correlations}
\end{figure}

\begin{figure}[ht!]
  \centering
  \includegraphics[width=.49\textwidth]{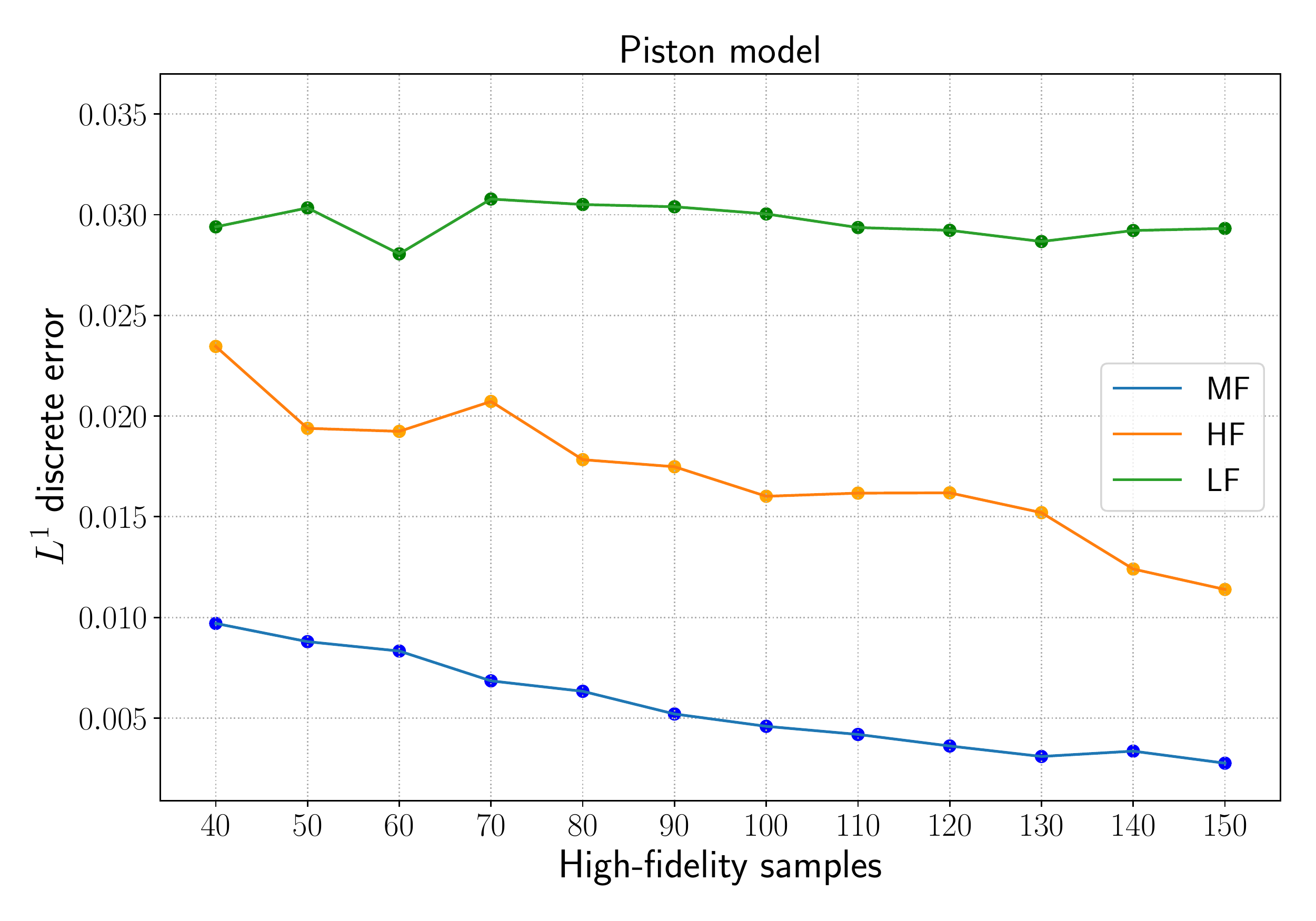}\hfill
  \includegraphics[width=.49\textwidth]{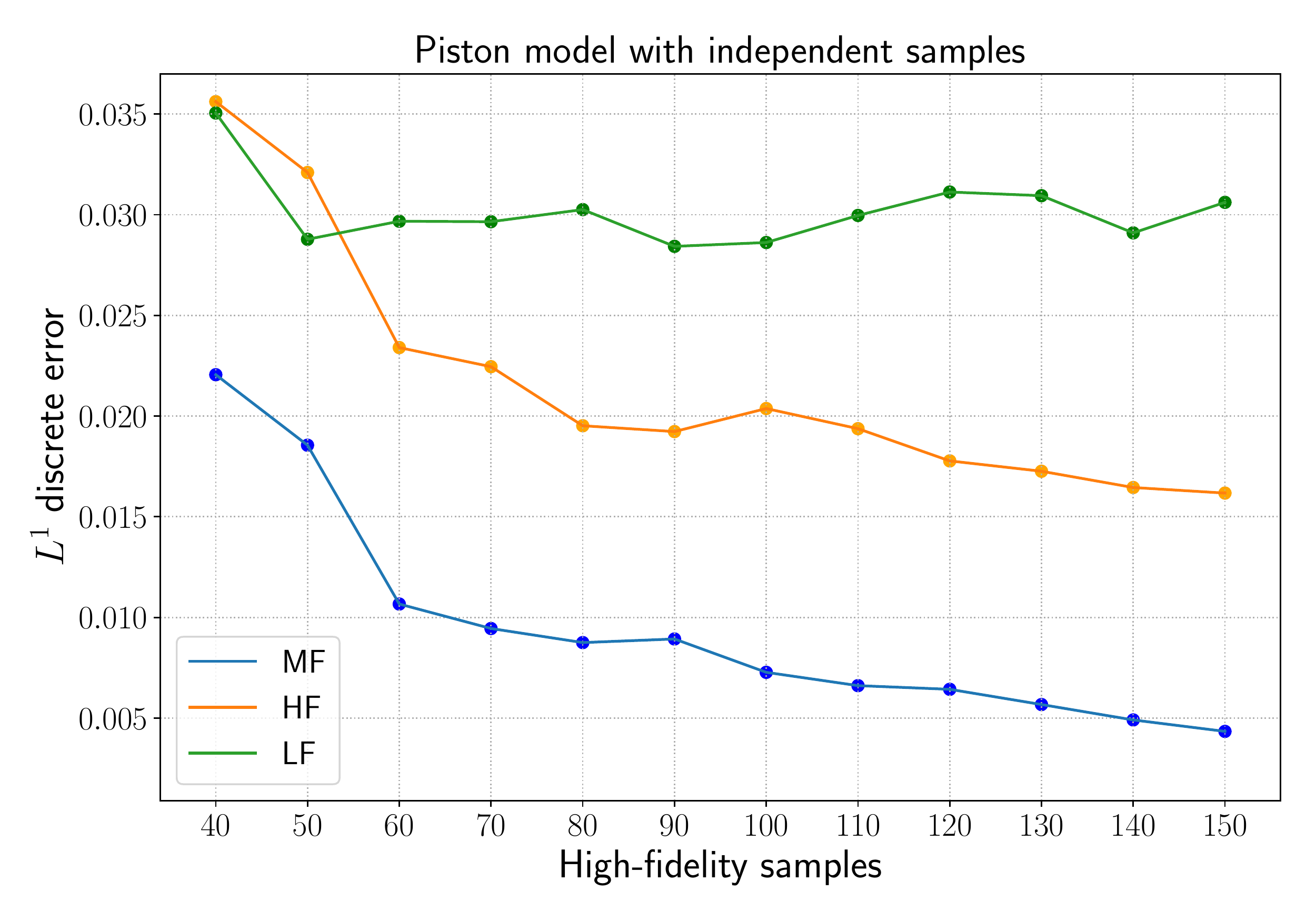}\\
  \caption{$L^{1}$ discrete error of the posterior of the multi-fidelity (MF),
    high-fidelity (HF) and low-fidelity (LF) models against the number of
    high-fidelity samples used to find the active subspace and build the
    Gaussian process regressions of the MF, HF, LF models. The 10000 test
    samples are distributed with Latin-Hypercube sampling. Left: identical
    high-fidelity samples for AS and NARGP. Right: independent high-fidelity
    samples for AS and NARGP.}
  \label{fig:piston}
\end{figure}

It is qualitatively evident from the sufficient summary plot in the
left panel of Figure~\ref{fig:piston_correlations} that a
one-dimensional active subspace is enough to explain with a fairy good
accuracy the dependence of the output from the $7$-dimensional
inputs. This statement could be supported looking at the ordered
eigenvalues of the correlation matrix of the gradients, which would
show a spectral gap between the first and the second eigenvalues.
We can also see in the right panel of
Figure~\ref{fig:piston_correlations} the correlations
scatter plot among the two different fidelity levels of the NARGP
model. The low-fidelity GP regression built on the dataset $\{(x^{L}_{i},
y^{L}_{i})\}_{i=1}^{N_L}$ performs already a good regression without
many outliers in the predictions evaluated at the test samples.

Figure~\ref{fig:piston} shows the errors of the MF models built with
different procedures, as described in Section~\ref{sec:mfas}. It can be seen that
using independent samples for the active subspace evaluation does not improve
the predictions obtained. Since in the right panel one third of the
high-fidelity samples are used to identify the active subspace, to
evaluate the differences between the two algorithms we have to compare
the $150$ samples on the right with the $100$ samples on the left, for
example. We can clearly notice that the two approches perform
almost the same.

\paragraph{SEIR model for Ebola}
The SEIR model for the spread of Ebola depends on $8$ parameters and the output
of interest is the basic reproduction number $R_0$. A complete AS analysis was
made in~\cite{diaz2018modified}, while a kernel-based active subspaces
comparison can be found in~\cite{romor2020kas}. The formulation is the
following:
\begin{equation}
  \label{eq:Ebola}
  R_0 =\frac{\beta_1 +\frac{\beta_2\rho_1 \gamma_1}{\omega} +
    \frac{\beta_3}{\gamma_2} \psi}{\gamma_1+ \psi},
\end{equation}
where the parameters range are taken from~\cite{diaz2018modified}.

\begin{figure}[ht!]
  \centering
  \includegraphics[trim=0 10 0 10, clip, width=.48\textwidth]{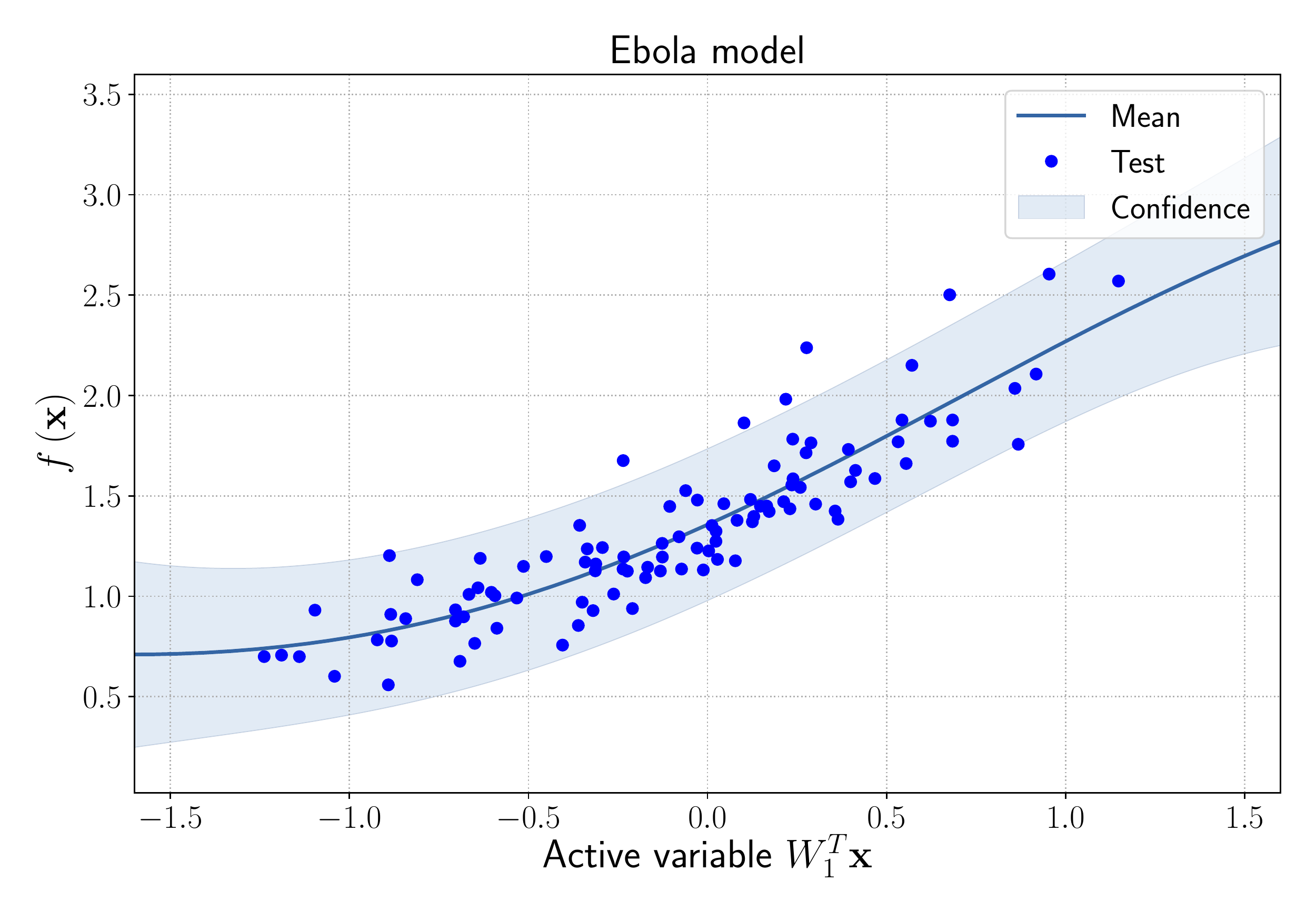}\hfill
  \includegraphics[trim=0 0 0 10, clip, width=.49\textwidth]{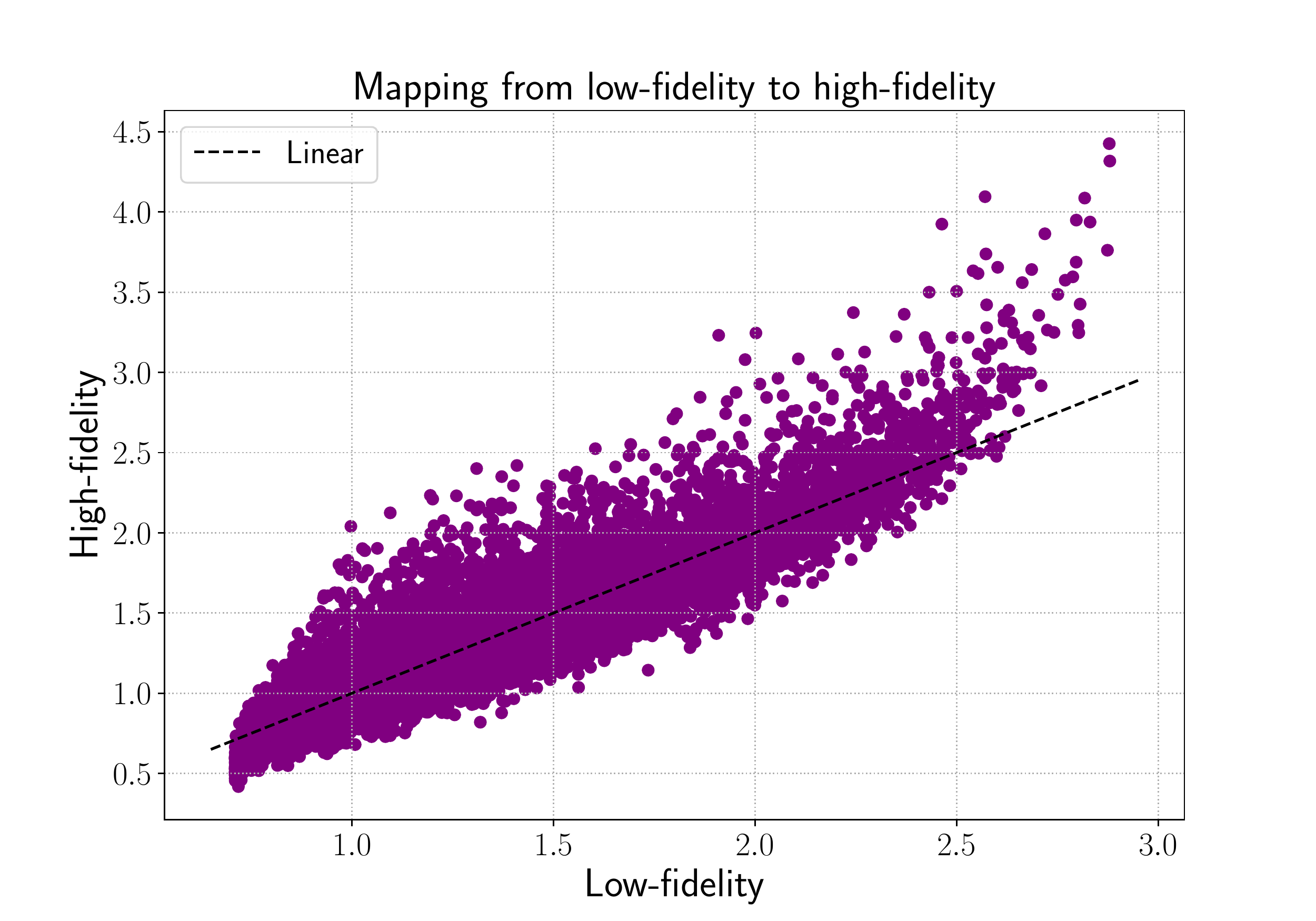}\\
  \caption{Left: sufficient summary plot of the Ebola model, 150
    samples were used to build the AS surrogate model shown. Right:
    correlation among the low-fidelity level and the high-fidelity
    level of the multi-fidelity model.}
  \label{fig:Ebola_correlations}
\end{figure}

\begin{figure}[ht!]
  \centering
  \includegraphics[width=.49\textwidth]{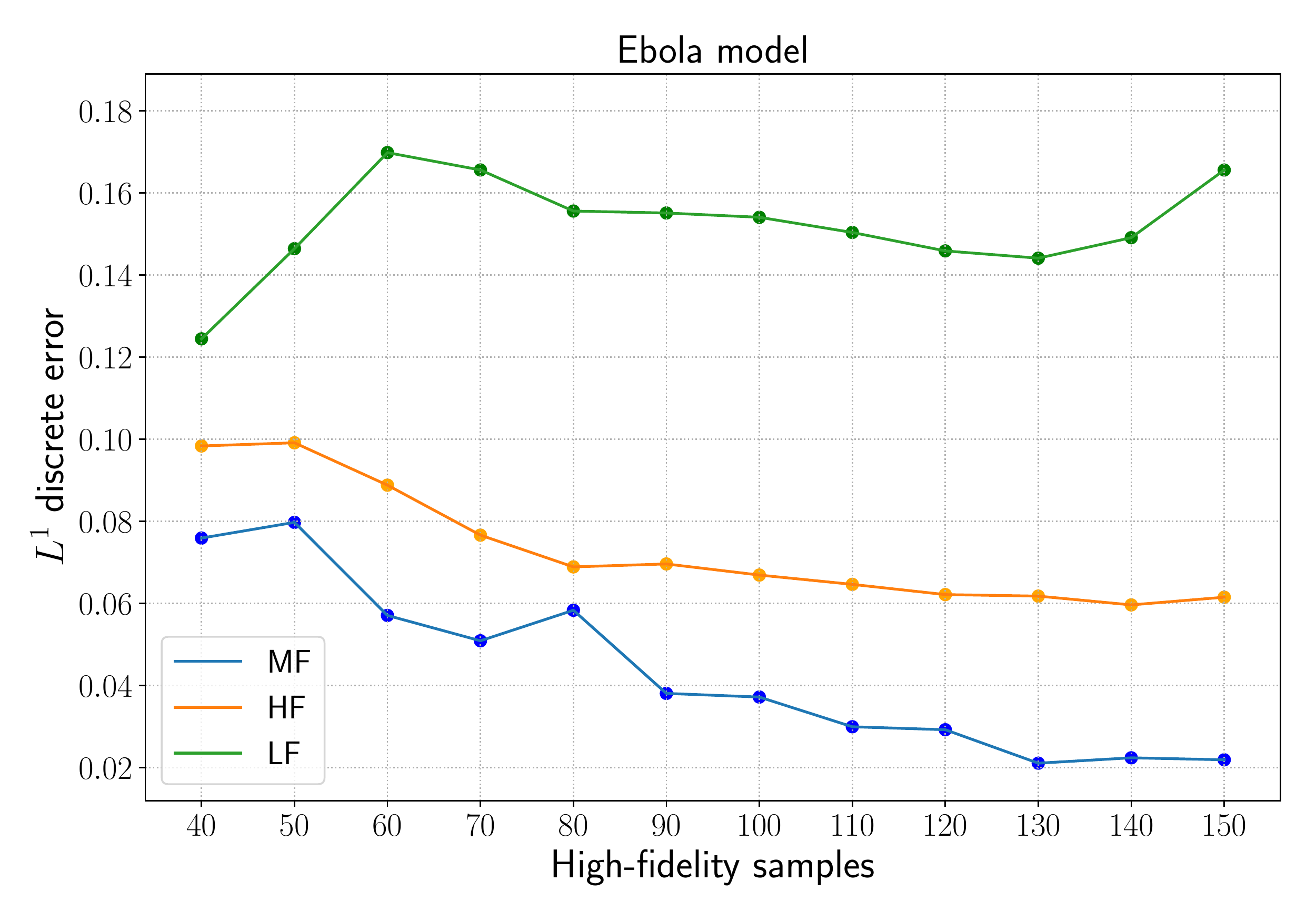}\hfill
  \includegraphics[width=.49\textwidth]{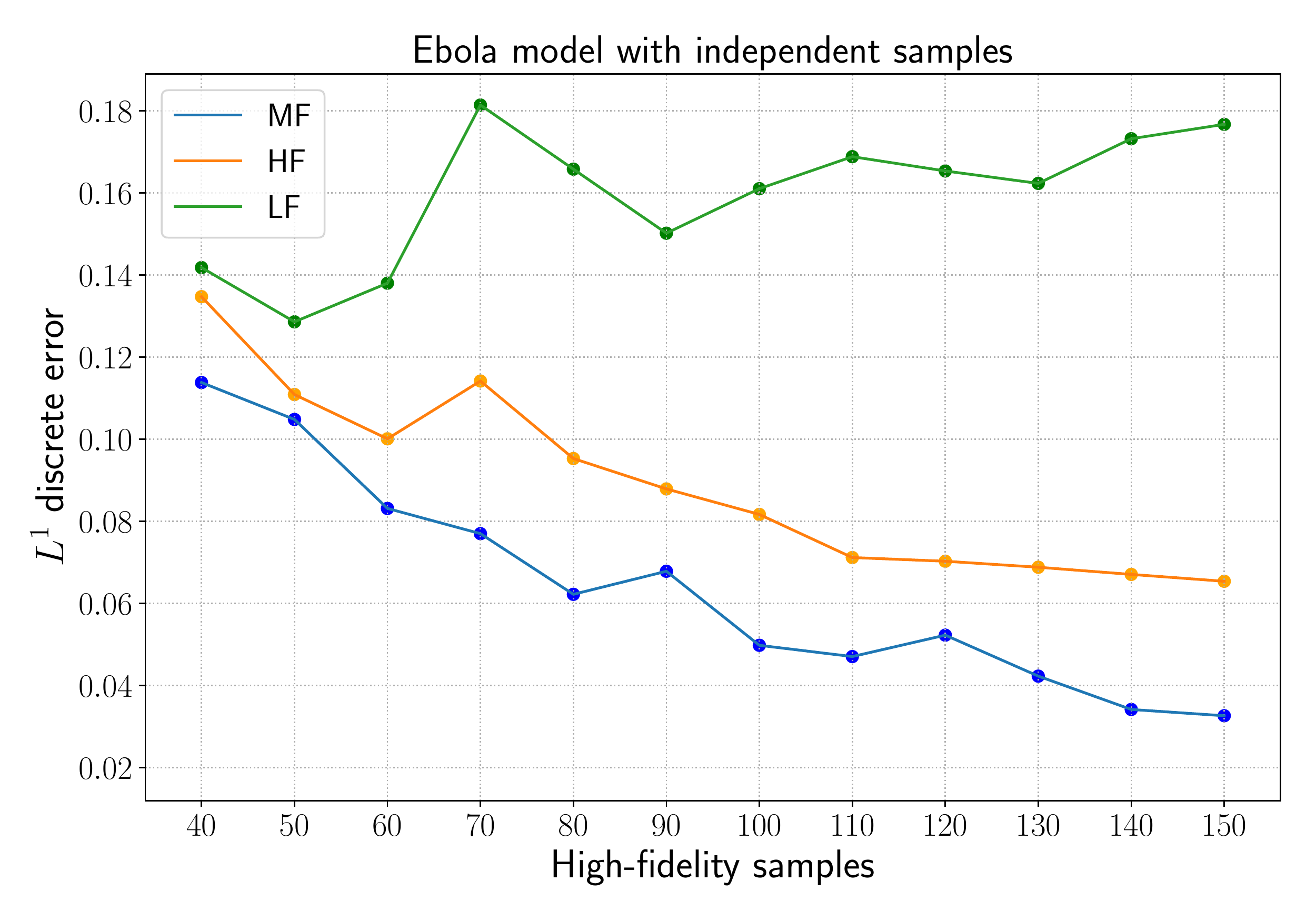}\\
  \caption{$L^{1}$ discrete error of the posterior of the multi-fidelity (MF),
    high-fidelity (HF) and low-fidelity (LF) models against the number of
    high-fidelity samples used to find the active subspace and build the
    Gaussian process regressions of the MF, HF, LF models. The 10000 test
    samples are distributed with Latin-Hypercube sampling. Left: identical
    high-fidelity samples for AS and NARGP. Right: independent high-fidelity
    samples for AS and NARGP.}
  \label{fig:Ebola}
\end{figure}

Differently from the previous test case, the one-dimensional response
function in the left panel of Figure~\ref{fig:Ebola_correlations} does
not explain well the model: in this case kernel-based active subspaces
could be employed to reach a better expressiveness of the surrogate
model~\cite{romor2020kas}. Even the scatter plot in the right panel
of Figure~\ref{fig:Ebola_correlations}, which shows the correlations between the
low-fidelity and high-fidelity levels of the NARGP model, exhibits a worse
accuracy in the low-fidelity level with respect to the previous test case. These
results are quantified in Figure~\ref{fig:Ebola} with the $L^1$
discrete error for the different fidelities models which are one order
of magnitude higher than the piston model, see Figure~\ref{fig:piston}.

From a comparison between the HF and MF models in Figure~\ref{fig:Ebola} and the
respective models in Figure~\ref{fig:piston} it can be seen that the nonlinear
autoregressive fidelity fusion approach learns relatively worse correlations of
the low-/high-fidelity levels of the NARGP Ebola model with respect to the
piston model.

For both the test cases the multi-fidelity regression approach with
active subspaces results in better performance with a consistent
reduction of the $L^1$ discrete error over a test dataset of $10000$ points.

\section{Conclusions}
\label{sec:conclusions}

In this work we propose a nonlinear multi-fidelity approach for the
approximation of scalar function with low intrinsic dimensionality. Such
dimension is identified by searching for the existence of an active subspace for
the function of interest. With a regression along the active variable we build a
low-fidelity model over the full parameter space which is fast to evaluate and
does not need any new simulations. We just extract new information from the
high-fidelity data we already have. This multi-fidelity approach results in a
decreased regression error and improved approximation capabilities over all the
parameter space.

We apply the multi-fidelity with AS method to two different benchmark problems
involving high dimensional scalar functions with an active subspace. We achieve
promising results with a reduction of the $L^1$ discrete error around $60$--$70\%$
with respect to the high-fidelity regression in one case (piston) and around
$20$--$66\%$ in the other one (Ebola), depending on the number of high-fidelity
samples used.

Further investigation will involve the use of more active subspaces based
fidelities, such as the kernel active subspaces~\cite{romor2020kas}. This could
also greatly improve data-driven non-intrusive reduced order
methods~\cite{morhandbook2019, rozza2018advances, salmoiraghi2016advances,
tezzele2018ecmi} through modal coefficients reconstruction and prediction for
parametric problems. We also mention the possible application to shape
optimization problems for the evaluation of both the target function and the
constraints.

\section*{Acknowledgements}
This work was partially supported by an industrial Ph.D. grant sponsored by
Fincantieri S.p.A. (IRONTH Project), by MIUR (Italian ministry for university
and research) through FARE-X-AROMA-CFD project, and partially funded by European
Union Funding for Research and Innovation --- Horizon 2020 Program --- in the
framework of European Research Council Executive Agency: H2020 ERC CoG 2015
AROMA-CFD project 681447 ``Advanced Reduced Order Methods with Applications in
Computational Fluid Dynamics'' P.I. Professor Gianluigi Rozza.

\bibliographystyle{abbrv}
\bibliography{pamm_bib.bib}

% \begin{thebibliography}{10} \end{thebibliography}

\end{document}